\documentclass{statsoc}
\usepackage{natbib}
\usepackage{epsfig}
\newtheorem{theorem}{Theorem}

\newtheorem{definition}{Definition}

\title[Functional Responses]{Nonlinear functional models for functional responses in reproducing kernel Hilbert spaces}
\author[H. Lian]{Heng Lian}
\runauthor{H. Lian}
\address{Brown University, Providence, USA.}
\coaddress{Heng Lian, Division of Applied Mathematics, 182 George St, Brown University, Providence, RI 02912, USA.}
\email{Heng\_Lian@brown.edu}

\begin{document}
\begin{abstract}
An extension of reproducing kernel Hilbert space (RKHS) theory provides a new framework for modeling functional regression models with functional responses. The approach only presumes a general nonlinear regression structure as opposed to previously studied linear regression models. Generalized cross-validation (GCV) is proposed for automatic smoothing parameter estimation. The new RKHS estimate is applied to both simulated and real data as illustrations. 
\end{abstract}
\section{Introduction}
In many experiments, functional data appear as the basic unit of observations. As a natural extension of the multivariate data analysis, functional data analysis provides valuable insights into these problems. Compared with the discrete multivariate analysis, functional analysis takes into account the smoothness of the high dimensional covariates, and often suggests new approaches to the problems that have not been discovered before. Even for nonfunctional data, the functional approach can often offer new perspectives on the old problem. 

The literature contains an impressive range of functional analysis tools for various problems including exploratory functional principal component analysis, canonical correlation analysis, classification and regression. Two major approaches exist. The more traditional approach, carefully documented in the monograph \citet{fda}, typically starts by representing functional data by an expansion with respect to a certain basis, and subsequent inferences are carried out on the coefficients. The most commonly utilized basis include B-spline basis for nonperiodic data and Fourier basis for periodic data. Another line of work by the French school, taking a  nonparametric point of view, extends the traditional nonparametric techniques, most notably the kernel estimate, to the functional case. Some theoretical results are also obtained as a generalization of the convergence properties of the classical kernel estimate. 

We are concerned with the regression problem in this paper. In general, the regression problem takes the form
\begin{equation}\label{model}
y_i=F(x_i)+\epsilon_i
\end{equation}
In traditional nonparametric inference, with $x_i$ and $y_i$ both being scalars, there exist a large number of methods including kernel and locally linear estimation. In functional data analysis, one or more of the components, $x_i, y_i$, and $\epsilon_i$, are functions defined on some interval, here assumed to be the interval $[0,1]$. Inferences are focused on the estimation of the structural component $F(x)$, with the residual $\epsilon$ modeling the noise or in general the component of observations not captured by the model. 

In the case of scalar responses $y$, at least two nonparametric approaches have appeared in the literature. The first method uses a simple kernel regression estimate
\[
\hat{F}(x)=\frac{\sum_i k(d(x_i,x))y_i}{\sum_i k(d(x_i,x))}
\]
where $d$ is a semi-metric on the space of functions. The second method is to use the RKHS framework, assuming the real valued function $F$ is an element of the RKHS $H$, and the estimator is obtained as the minimizer of the regularized loss
\begin{equation}\label{rkhsscalar}
\sum_i(y_i-F(x_i))^2+\lambda ||F||_H
\end{equation}
The construction of the Hilbert space in this scalar response case involves no extra technical difficulties compared with the classical multivariate case as long as we have a metric on the space of functions $x$, and then the kernel can be constructed with $K(x_1,x_2)=k(d(x_1,x_2))$ for any positive definite function $k$. 
The representer theorem for RKHS implies that the solution to (\ref{rkhsscalar}) has the form \[ 
F(x)=\sum_i \alpha_ik(x_i,x)
\]
with real coefficients $\alpha_i$. This representation can be plugged back into (\ref{rkhsscalar}) and solve for the coefficients. Note that in both of these two nonparametric approaches, the same word ``kernel" is used for different concepts, the exact meaning of the word should be clear from the context.

In the case of functional responses $y$, the classical parametric inferences assume that $F$ is  linear in $x$. More explicitly, the pointwise model assumes that 
\[
y(t)=\alpha(t)+\beta(t)x(t)+\epsilon(t)
\]
while the integral model specifies that
\begin{equation}\label{linest}
y(t)=\alpha(t)+\int_0^1 \beta(s,t)x(s)\,ds+\epsilon(t)
\end{equation}
In contrast to traditional linear regression models, now the coefficient $\beta$ is a function on $[0,1]$, or a bivariate function on $[0,1]\times [0,1]$. To estimate the coefficient $\beta$, again a basis is chosen to represent the functions involved and the problem is reduced to a multiple linear regression model. To our knowledge, nonparametric approaches to functional analysis with functional responses has not been studied before and we will embark on this task in the current paper.

We consider the functional response model (\ref{model}) in which the structural component $F$ is a mapping from some space of functions to another function space. We assume that $y(\cdot)$ belongs to a Hilbert space $H$. Although it is not necessary to assume that $x(\cdot)$ belong to the same space, or even that there is an inner product associated with it, it will be convenient to require that $x$ is in the same space as $y$, as we will assume in the following.

In this paper, we will develop an estimation procedure for functional response models within the framework of RKHS. For nonlinear models, the parametric approach above does not give satisfactory results. Our goal is to show that within the RKHS framework, a simple estimate can be developed which reduces to linear regression computations very much like those derived in the parametric approach. Our work is motivated by the work of \citet{preda}. Unlike the case of scalar response models treated in \citet{preda}, the extension we need is more complicated, as we will discuss in Section 2. There we show how we should extend the notion and the construction of a RKHS in this new setting and also prove the corresponding representer theorem. In section 3, we present our new nonlinear model within the framework developed and also comment on the computations involved. Simulation studies and application to the well-known weather data are carried out in section 4. These results show clear advantage of our model compared to the parametric linear regression model in nonlinear contexts. A kernel estimate similar to that of \citet{ferraty2003} and \citet{ferraty2004} are also constructed for comparison. We conclude the paper in Section 5. Some technical details are deferred to the appendix.

\section{Reproducing kernel Hilbert spaces}
Following \citet{wahba}, a (real) RKHS $H$ is a Hilbert space of real-valued functions defined on, say, the interval $[0,1]$, in which the point evaluation operator $L_t: H\rightarrow R, L_t(f)=f(t)$ is continuous. By Riesz representation theorem, this definition implies the existence of a bivariate function $K(s,t)$ such that
\begin{eqnarray}\label{rep}
&&K(s,\cdot)\in H, \mbox{ for all } s\in [0,1]\nonumber\\
&&\mbox{(reproducing property) for every } f\in H \mbox{ and } t\in [0,1],\;  \langle K(t,\cdot), f\rangle_H=f(t)
\end{eqnarray}
The definition of a RKHS can actually start from a positive definite bivariate function $K(s,t)$ and RKHS is constructed as the completion of the linear span of $\{K(s,\cdot),s\in [0,1]\}$with inner product defined by $\langle K(s,\cdot),K(t,\cdot)\rangle_H=K(s,t)$.

In the regression model (\ref{model}) with functional response, we are dealing with the functional $F$ taking values in the Hilbert space $H$. So the RKHS we construct should be a subset of $\{F: H\rightarrow H \}$. To define the RKHS in this case, we follow the same procedure as in \citet{wahba}. 
\begin{definition}
A (functional) RKHS $\mathcal{H}$ is a subset of $\{F: H\rightarrow H \}$. It is a Hilbert space, with inner product $\langle \cdot,\cdot\rangle_\mathcal{H}$, in which the point evaluation operator is a bounded linear operator, i.e., $L_x: F\rightarrow F(x)$ is a bounded operator from $\mathcal{H}$ to $H$ for any $x \in H$.
\end{definition}

The definition above is not useful for constructing a RKHS. For our purpose, we need to explicitly define the kernel associated with the space. The continuity of $F\rightarrow F(x)$ for any $x\in H$ is equivalent to the continuity of the mapping $F\rightarrow \langle F(x),g\rangle_H$ for any $x\in H$ and $g\in H$. By Riesz representation theorem applied to the Hilbert space $\mathcal{H}$, there exists an element $K_x^g$ in $\mathcal{H}$ such that
\begin{equation}\label{frep}
\langle K_x^g,F\rangle_\mathcal{H}=\langle F(x),g\rangle_H
\end{equation}
From the above, the mapping $g\rightarrow K_x^g$ is linear. By the boundedness of the point evaluation operator, we also have
\[
\langle K_x^g,F\rangle_\mathcal{H}=\langle F(x),g\rangle_H\le C||F||_\mathcal{H}||g||_H
\]
which implies that the mapping $g\rightarrow K_x^g$  is also bounded. And so $g\rightarrow K_x^g(y)$ is a bounded linear operator for any $y\in H$, which we can denote by $K(x,y)$, i.e, $K(x,y)g:=K_x^g(y)$, and we call $K(\cdot,\cdot)$ the reproducing kernel associated with $\mathcal{H}$. Note that in this case, the reproducing property is defined by $\langle K(x,\cdot)g,F\rangle_\mathcal{H}=\langle F(x),g\rangle_H$ for any $x \in H, g\in H$, this is just a rewriting of (\ref{frep}).

Two properties of the reproducing kernel are immediate. (a) $K(x,y)=K(y,x)^{T}$, where the superscript $T$ denotes the adjoint operator. This is a simple consequence of the following sequence of identities making use of (\ref{frep}): $\langle K(x,y)g,f\rangle_H=\langle K_x^g(y),f\rangle_H=\langle K_x^g,K_y^f\rangle_\mathcal{H}=\langle K_y^f(x),g\rangle_H=\langle g,K(y,x)f\rangle_H $. (b) $\sum_{i=1}^n\sum_{j=1}^n \langle K(x_i,x_j)f_i,f_j\rangle_H\ge 0$, which follows from $\sum_{i,j}\langle K(x_i,x_j)f_i,f_j\rangle_H=\sum_{i,j}\langle K^{f_i}_{x_i},K^{f_j}_{x_j}\rangle_\mathcal{H}=||\sum_iK^{f_i}_{x_i}||_\mathcal{H}$.

From the above discussions, we have the following definition of a positive definite kernel.
\begin{definition}
A (functional) nonnegative definite kernel is a bivariate mapping on H$\times$H, taking values in $L(H)$, the space of bounded linear operators from $H$ to itself, such that
\[K(x,y)=K(y,x)^T\] and
\[\sum_{i=1}^n\sum_{j=1}^n \langle K(x_i,x_j)f_i,f_j\rangle_H\ge 0\]
where $ \{x_i\} $ and $ \{f_i\} $ are any two sequences in $H$. If the double sum above is strictly positive whenever $\{x_i\}$ are distinct and $\{f_i\}$ are not all zero, $K$ is a positive definite kernel.
\end{definition}

Given a positive definite kernel defined as above, we can construct a unique RKHS of functions on $H$ taking values also in $H$, with $K$ as its reproducing kernel. The proof of this statement follows exactly the same lines as for the real-valued case.

\section{Models for functional data}
We consider the inference and prediction problem for model (\ref{model}). Given the observed functional covariates and responses $\{(x_1,y_1), \ldots, (x_n,y_n)\}$, a general approach to estimate $F$ is to solve the following minimization problem:
\begin{equation}\label{hyperopt}
\min_{F\in \mathcal{H}}\sum_{i=1}^n||y_i-F(x_i)||_2^2+\lambda||F||_\mathcal{H}
\end{equation}
where a penalty term in $\mathcal{H}$-norm with smoothing parameter $\lambda>0$ is added to the empirical risk as is usually done in the smoothing spline regression. We use the simplest loss function $||x||_2^2=\int_0^1x^2(t)\,dt$ here, although other types of loss can certainly be considered.

The optimization problem above optimizes over the infinite-dimensional space $\mathcal{H}$ which is not feasible in general. Fortunately, the representer theorem below reduces this problem to finite dimensions (if you consider $\mathcal{H}$ as a vector space with elements in $H$ acting as ``constants"). The proof of the following is deferred to the appendix.
\begin{theorem}
Given the observations $\{(x_i,y_i)\}_{i=1}^n$, the solution to (\ref{hyperopt}) has the following representation
\begin{equation}\label{hyperrep}
\hat{F}=\sum_{i=1}^n K(x_i,\cdot)\alpha_i
\end{equation}
with functional coefficients $\alpha_i\in H$.
\end{theorem}

Two difficulties arise at this stage. First, the construction of $K$ in general is difficult and a search of the literature does not seem to provide us with any clues about how to construct a positive definite kernel in general. Second, even if the kernel is constructed, it is not clear how to solve (\ref{hyperopt}) after we plug in the representation (\ref{hyperrep}). Due to these difficulties, we are forced to choose the simplest functional kernel $K(x,y)=a(d(x,y))I$, where $a(\cdot)$ is a real positive definite function, $I\in L(H)$ is the identity operator, and $d$ is a metric on $H$. We choose the simplest metric $d(x,y)=||x-y||_2$ which is also used in \citet{preda}. It is unfortunate that we will not be able to provide more complicated examples of the kernel, but this estimate works reasonably well in all our experiments. It is clear that $K$ defined in this way is nonnegative definite, but to prove that it is positive definite requires a little more extra work. We state this result formally as a theorem and its proof is to be found in the appendix.
\begin{theorem}
The functional kernel $K(x,y)=a(||x-y||_2)I$ is positive definite if $a(\cdot)$ is a (real) positive definite function.
\end{theorem}

After this dramatic simplification, we are able to solve problem (\ref{hyperopt}). Let $a_{ij}=a(||x_i-x_j||_2)$ and make use of the representer theorem 1, we arrive at the following optimization problem
\begin{equation}\label{opt}
\min_{\alpha_i}\sum_{i=1}^n||y_i-\sum_j a_{ij}\alpha_j||_2^2+\lambda\sum_{i,j}a_{i,j}\langle\alpha_i,\alpha_j\rangle_H
\end{equation}

From this point on, there are definitely more than one way to preceed. For example, one can represent each $\alpha_i$ by expansion with respect to a chosen basis. We choose instead to again compute (\ref{opt}) in the RKHS framework by assuming $H$ is itself a (real) RKHS with reproducing kernel $k$. The loss term is first discretized, assuming the observations are made on a regular grid $\{t_1,\ldots,t_T\} $ on $[0,1]$, then  another application of the representer theorem in the real-value case stated in Theorem 3 below reveals the solution to be 
\begin{equation}\label{discreterep}
\alpha_i=\sum_{l=1}^T b^i_l k(t_l,\cdot)
\end{equation}
and we only need to compute the coefficients $\{b^i_l\}$. Formally, we have
\begin{theorem}
Consider the raw data version of  (\ref{opt})
\begin{equation}\label{discreteopt}
\min_{\alpha_i}\sum_{i=1}^n\sum_{l=1}^T[y_i(t_l)-\sum_j a_{ij}\alpha_j(t_l)]^2+\lambda\sum_{i,j}a_{i,j}\langle\alpha_i,\alpha_j\rangle_H
\end{equation}
The solution to (\ref{discreteopt}) is of the form  (\ref{discreterep}).
\end{theorem}

The readers are referred to the appendix for detailed formula involved in the computation.

The model (\ref{discreteopt}) looks similar to the classical smoothing spline estimation with two differences. First, in the first term of (\ref{discreteopt}), instead of trying to smooth $y$ by a single function, it tries to approximate each observation $y_i$ as a combination of a common set of functions $\{\alpha_j\}_{j=1}^n$. The coefficients $a_{ij}$ reflects the similarity of the covariates $x_i$ and $x_j$. Second, an unconventional penalty term $\sum_{i,j}a_{ij}\langle\alpha_i,\alpha_j\rangle_H$ appears, which involves cross-over terms $\langle\alpha_i,\alpha_j\rangle_H$ with $i\neq j$.  The loss term in (\ref{discreteopt}) seems natural and one could probably come up with this term without going through all the trouble of using functional RKHS and the derivation above. If this is the case, one would probably use a penalty term such as $\sum_{i}a_{ii}||\alpha_i||_H$ for regularization purposes. We call the solution to the problem (\ref{discreteopt}) with this simpler penalty the modified RKHS estimate, as opposed to the RKHS estimate which solves the original problem (\ref{discreteopt}). 

In our implementation, we use the Gaussian kernel for both RKHSs, so $K(x,y)=exp\{-||x-y||_2^2/2\sigma^2\}I$ and $k(s,t)=exp\{-(s-t)^2/2\sigma'^2\}$. Thus there are at least three parameters, $\sigma,\sigma',\lambda$, that need to be specified. For the width parameters in the kernels, we simply choose $\sigma$ to be the mean of all $||x_i-x_j||_2, i,j=1,\ldots,n$, and $\sigma'$ is similarly chosen to be the mean of $||t_i-t_j||_2, i,j=1,\ldots,T$. These choices are of course suboptimal but produce good results in our experiences, so we do not bother to search over these parameter spaces. The choice of the smoothing parameter is arguably more important. Generalized cross-validation (GCV), which has been extensively studied in \citet{wahba}, can be used for the selection of $\lambda$. Given a grid of values for $\lambda$ specified beforehand, GCV approximates the true error by 
\[
V(\lambda)=\frac{1}{n}||(I-A(\lambda))y||^2/[\frac{1}{n}Tr(I-A(\lambda))]^2
\]
where $A(\lambda)$ is the ``influence matrix" defined in Appendix B, and the final $\lambda$ is chosen to be the one that minimizes $V(\lambda)$.

A nonparametric kernel estimate is studied in \citet{ferraty2003}, \citet{ferraty2004}, and \citet{ferraty2007}, which is simply defined as 
\begin{equation}\label{kernelest}
\hat{F}(x)=\frac{\sum_i k(||x_i-x||)y_i}{\sum_ik(||x_i-x||)}
\end{equation}
In those papers, the authors only studied the kernel estimate for the model with scalar responses, but this estimate can obviously be used when the dependent variable is a curve. It also seems natural that we should smooth the response $y_i$ before plugging it into (\ref{kernelest}) if the observations are noisy.  In the next section, this kernel estimate will be compared with our RKHS estimates and the linear parametric estimate.

\section{Examples}
\subsection{Simulation study}
One of the main goals of using RKHS estimate is to deal with nonlinearity in the data. In this simulation study, we compare four estimates. The RKHS estimate (\ref{hyperrep}), its simpler version with modified penalty, the linear regression estimate (\ref{linest}) and the kernel estimate (\ref{kernelest}). The simulated data are generated by the following models:
\begin{enumerate}
\item $y(t)=\int_0^1|t-s|x(s)\,ds$
\item $y(t)=\int_0^1|t-s|x^2(s)\,ds$
\item $y(t)=\sin(2\pi t)x(t),\;t\in[0,1]$
\item $y(t)=\cos(\pi t)|x(t)|,\;t\in[0,1]$
\end{enumerate}
Models (a) and (c) are linear in $x$ while (b) and (d) are nonlinear. $x(\cdot)$ is generated as a standard Brownian motion with random starting position uniform in $[0,5]$, and $y(t)$ is computed using expressions (a)-(d). An equispaced grid of 50 points on $[0,1]$ is used. The raw observations for the dependent variable are the values of $y$ on the grid points with i.i.d. standard normal noise added. The width parameters $\sigma$ and $\sigma'$ in RKHS and modified RKHS estimates are set to be the mean of the distances of the covariates from the training data as explained in the last section. After some experimentation, we use B-spline basis of order 4 with 10 equispaced knots in the linear regression estimate (\ref{linest}), with a penalty term involving the second partial derivatives of $\beta(s,t)$. The fitting of this linear model is performed using the fda package provided in R. A total of $n=30$ replicates are used in model fitting. Since this is a simulation study, we generate a separate set of 50 replicates as validation to choose the smoothing parameter $\lambda$ in all three models, so GCV is not used here. For the kernel estimate, we again used a Gaussian kernel, and the only parameter is the width parameter inside the kernel. Although it can be fixed as in the RKHS estimator, and it indeed produces good results, we instead search over this one-dimensional space using the same validation data that was used to choose the smoothing parameter in other models. In kernel estimates, we use the raw data as $y_i$ in (\ref{kernelest}). In simulation studies, we can also use $y(\cdot)$ in (\ref{kernelest}) before the standard normal noises are added, which we call the oracle kernel estimate for obvious reasons. In real applications, the performance of the kernel estimate should be worse than the oracle estimate if some kind of smoothing on the raw dependent data are use as a preprocessing step. The search over the width parameter put the kernel estimates in a favorable position compared to other estimates. Finally, after the parameters are estimated, we generate another 50 observations from model (a)-(d) as test data. These simulations are repeated 50 times for each model (a)-(d), and the mean of the mean square errors are reported in Table \ref{tabsim}.
\begin{table}
\caption{\label{tabsim}Results for the simulation study}
\fbox{\begin{tabular}{cccrcc}
&RKHS& Modified RKHS&Linear&Oracle Kernel&Kernel\\
\hline
Model (a)&1.000 &0.944&0.327&0.910 &2.773\\
Model (b)&1.000 &1.000&11.632&1.984 &2.035\\
Model (c)& 1.000&0.992&0.471&1.643 &1.954\\
Model (d)& 1.000&1.007&2.892&1.459 &1.735\\
\end{tabular}}
\end{table}

In Table \ref{tabsim}, the errors for the RKHS estimates are normalized to be $1$, and the errors of other estimates are shown as relative to the error for RKHS estimates. For linear models (a) and (c), the linear estimate clearly wins. For nonlinear models (b) and (d), the linear estimate perform badly compared to other estimates. The performance of the RKHS estimate and the modified RKHS estimate are almost identical. Although the kernel estimates do not perform as good as the RKHS estimates, it deserves further investigations due to its low computational complexity.
\begin{figure}
\centering
\makebox{\epsfig{figure=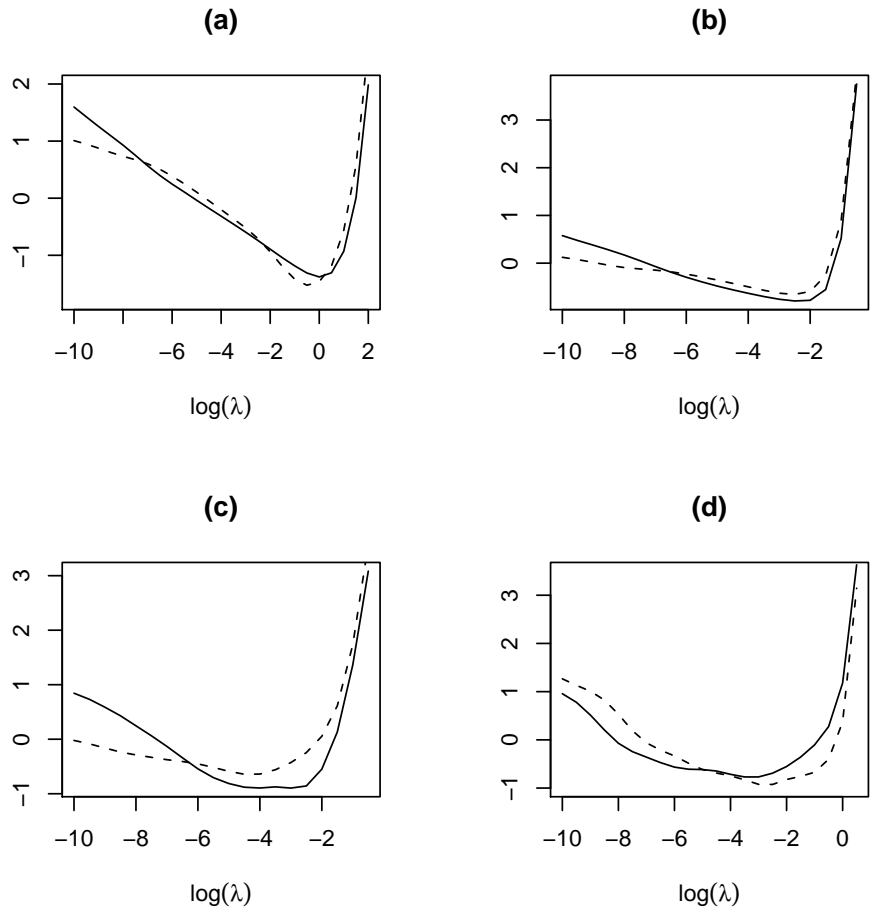,width=13cm}}
\caption{\label{gcv}Comparison of the GCV with the error computed from validation data. The solid curves are the GCV estimates, the dashed curves are the error computed from validation set. The curves are shifted and normalized to show the shape of the curves rather than its absolute magnitudes.}
\end{figure}

Next we study the performance of GCV for estimating the smoothing parameter in the RKHS estimate and compare it to the estimate obtained from the validation data, which act as the background truth in our simulation study. Figure \ref{gcv} show that for all four models (a)-(d), the GCV correctly identifies a good smoothing parameter to use in these cases.

\subsection{Application to the weather data}
\begin{figure}
\begin{minipage}[b]{.42\linewidth}
  \centering
  \centerline{\epsfig{figure=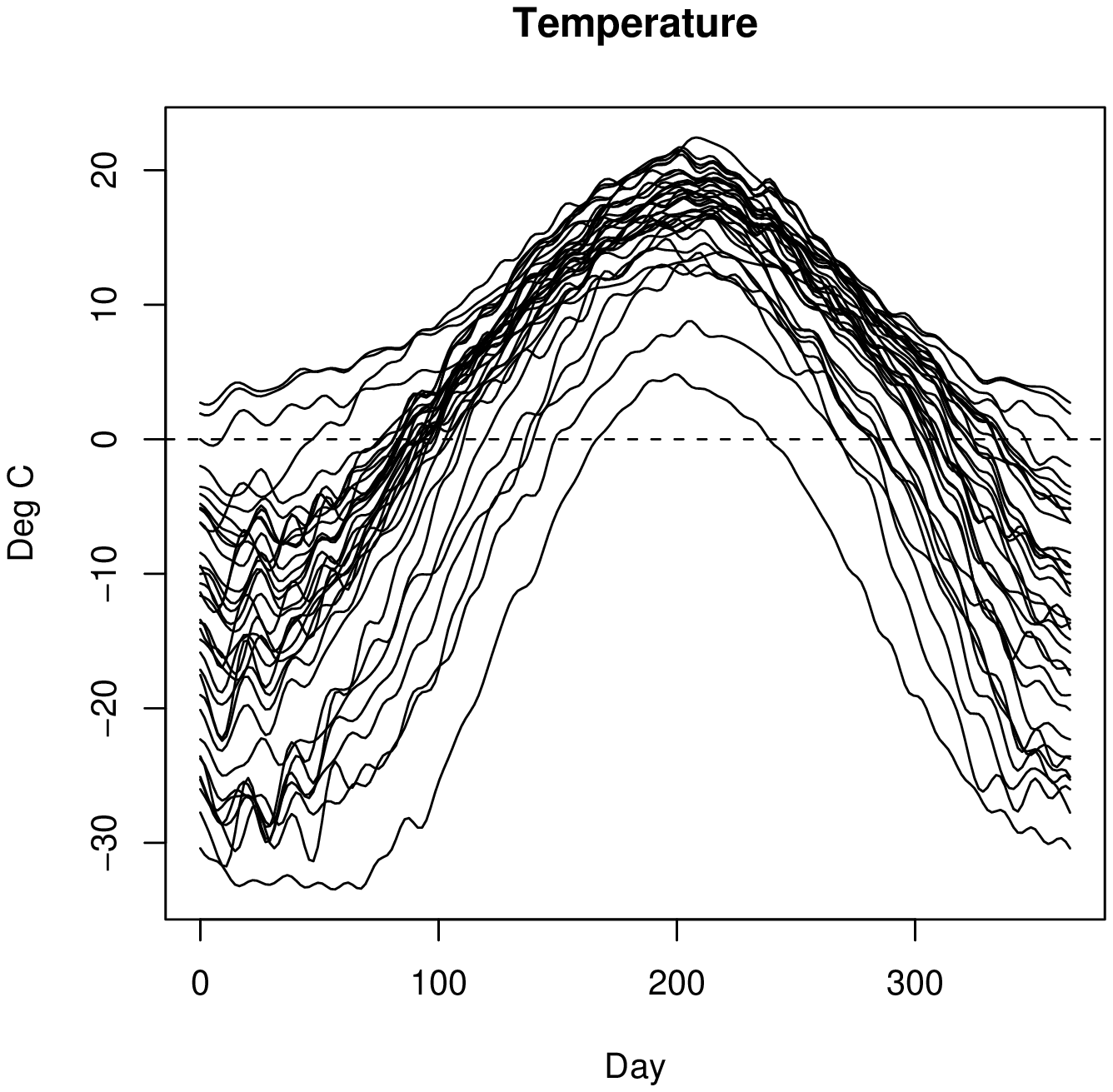,width=8.0cm}}
\end{minipage}
\hfill
\begin{minipage}[b]{0.42\linewidth}
  \centering
  \centerline{\epsfig{figure=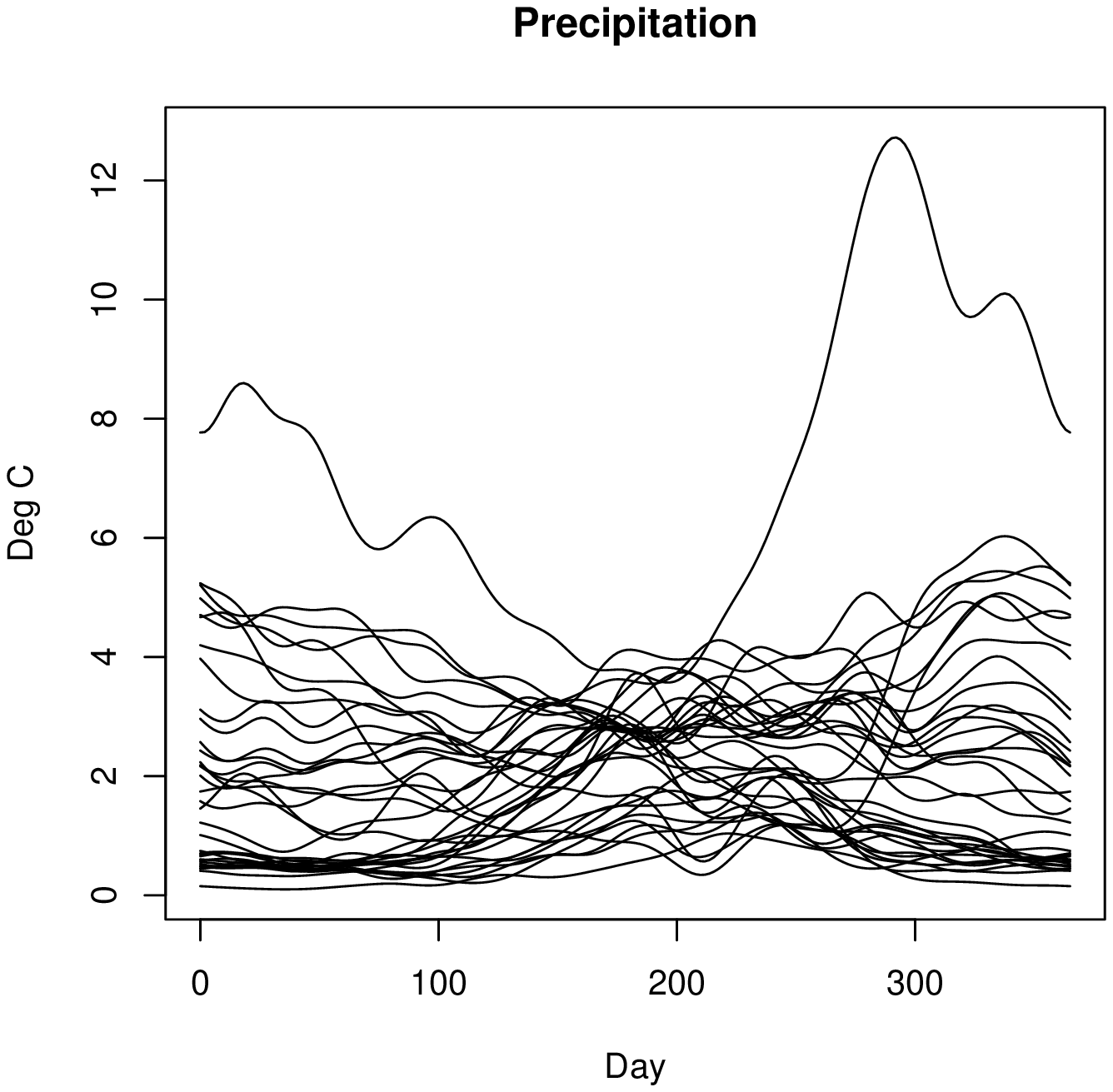,width=8.0cm}}
\end{minipage}
\caption{\label{weatherdata}Daily weather data for 35 Canadian stations, the curves plotted here result from using smoothing splines to fit the raw data.}
\end{figure}

\begin{figure}
\centering
\makebox{\includegraphics{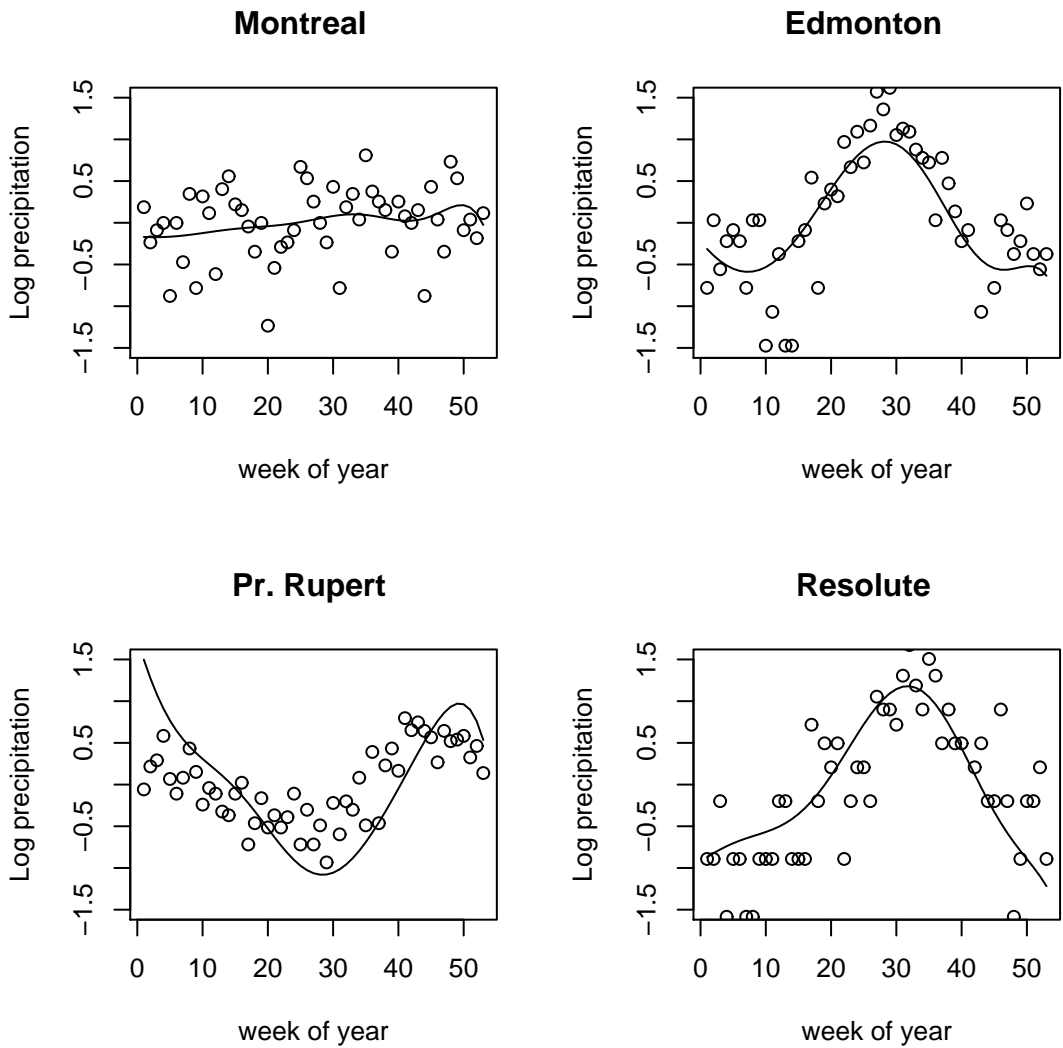}}
\caption{\label{weatherpredict}Raw data (points) and predictions (solid) of log precipitation for four weather stations.}
\end{figure}

The daily weather data consists of daily temperature and precipitation measurements recorded in 35 Canadian weather stations. These data are plotted in Figure \ref{weatherdata}. To save on computations, we subsample the data and use only the weekly measurements, so each observation consists of functional data observed on a equispaced grid of 53 points. We treat the temperature as the independent variable and the goal is to predict the corresponding precipitation curve given the  temperature measurements. As is previously done, we set the dependent variable to be the log tranformed precipitation measurements, and a small positive number is added to the values with $0$ precipitation recorded. The prediction of our RKHS estimate is shown in Figure \ref{weatherpredict} for four weather stations. Each plot is produced by fitting the model on the other 34 replicates, using GCV to choose the smoothing parameter, and then finally calculating the predicted precipitation based on the temperature curve. The figure shows reasonable prediction accuracy and can be compared to the results presented in Chapter 16 of \citet{fda}.

\section{Conclusions}
We proposed a new approach to fitting a nonlinear functional regression model for functional responses. The simulations we conducted demonstrated the clear advantage of the RKHS estimate over the linear regression model when the true model is nonlinear. The estimate is also better than the simplistic kernel estimates used in traditional nonparametric regression. 

The advantages of the RKHS estimates are tied with the additional computational costs. In problem (\ref{discreteopt}), we are optimizing over the same number of curves as there are the number of replicates which is computationally difficult when $n$ is large. Approximate solutions such as choosing a limited number of kernel basis centered around selected covariates may prove to be useful.

\appendix
\section*{Appendix A Proofs for Theorem 1-3}
Proof of Theorem 1. Denote by $\mathcal{H}_0$ the subspace of $\mathcal{H}$ spanned by the kernel centered at the observed covariates: $\mathcal{H}_0=\{\sum_{i=1}^nK(x_i,\cdot)\alpha_i, \alpha_i\in H\}$. Any $F\in \mathcal{H}$ can be written as $F=F_0+G$, where $F_0\in\mathcal{H}_0$  and $G\perp\mathcal{H}_0$ in $\mathcal{H}$. Then for any $j\in\{1,\ldots,n\}$, and any $h\in H$, $\langle G(x_j),h\rangle_H=\langle K(x_j,\cdot)h,G\rangle_\mathcal{H}=0$ by the reproducing property (\ref{frep}), this implies $G(x_j)=0$ by the arbitrariness of $h\in H$. Also, by the orthogonality of $G$ and $F_0$, $||F||_\mathcal{H}=||F_0||_\mathcal{H}+||G||_\mathcal{H}>||F_0||_\mathcal{H}$ if $G$ is nonzero. This implies that 
\[
\sum_{i=1}^n||y_i-F(x_i)||_2^2+\lambda||F||_\mathcal{H}>\sum_{i=1}^n||y_i-F_0(x_i)||_2^2+\lambda||F_0||_\mathcal{H}
\]
if $G$ is nonzero. So the minimizer $\hat{F}$ belong to $\mathcal{H}_0$

Proof of Theorem 2. We want to prove that $\sum_{i,j}\langle K(x_i,x_j)f_i,f_j\rangle_H\ge 0$.  Since $ \langle K(x_i,x_j)f_i,f_j\rangle_H =a_{ij}\langle f_i,f_j\rangle_H$, the nonnegativity  follows immediately from Schur's Lemma, which asserts that the Hadamard product of two nonnegative matrices is again a nonnegative matrix. 

If $\sum_{i,j}\langle K(x_i,x_j)f_i,f_j\rangle_H=0$, and $\{x_i\}$ are distinct. Let $A=\{a_{ij}\}, B=\{\langle f_i,f_j\rangle_H\}$. $A$ is positive definite and $B$ is nonnegative definite. We have the factorization $B=E^TE$. So $\sum_{i,j}\langle K(x_i,x_j)f_i,f_j\rangle_H=\sum_{i,j}A_{ij}B_{ij}=\sum_{i,j}A_{ij}(E^TE)_{ij}=\sum_{i,j,k}A_{ij}E_{ki}E_{kj}=\sum_k (\sum_{i,j}E_{ki}A_{ij}E_{kj})$. By the positive definiteness of $A$, we get $E=0$, which in turn implies $f_i=0$ for all $i$.

Proof of Theorem 3. This is similar to the classical proof with smoothing splines. We write $\alpha_i=\alpha_{i,0}+g_i$, with $\alpha_{i,0}\in span\{k(t_l,\cdot)\}$ and $g_i$ in its orthogonal complement. Similar to the proof of theorem 1, we only need to show that $\sum_{i,j}a_{ij}\langle\alpha_{i,0}+g_i,\alpha_{j,0}+g_j\rangle_H\ge\sum_{i,j}a_{ij}\langle\alpha_{i,0},\alpha_{j,0}\rangle_H$ with equality only if $g_i=0$ for all $i$. The proof of this statement is contained in the proof of Theorem 2.
\section*{Appendix B Computational details}
We detail the computations involved in solving the model (\ref{discreteopt}). The calculations are very similar to those used in linear regression model (\ref{linest}) using basis expansion, which is not surprising due to the representer theorem. Let $A=\{a_{ij}\}$ and $B$ be the $n$ by $T$ matrix containing the coefficients in the expansion (\ref{discreterep}): $B=\{b^i_l\}$. Also, denote the $n\times T$ matrix $\{y_i(t_l)\}$ by $Y$ and $T\times T$ matrix $\{k(t_i,t_j)\}$ by $K$. The objective function (\ref{discreteopt}) can be written in matrix form:
\[
Tr((Y-ABK)(Y-ABK)^T)+\lambda Tr(ABKB^T)
\]
Taking the derivative with respect to B, we want to solve
\[
(A^TA)B(KK^T)+\lambda ABK=AYK
\]
vectorizing the above equation gives us a system of linear equations
\[
[(KK^T)\otimes(A^TA)+\lambda K\otimes A]vec(B)=(K\otimes A)vec(Y)
\]
So we have the formula\[vec(B)=[K\otimes A+\lambda I]^{-1}vec(Y)\]
and the fitted values are $\hat{Y}=ABK$, the vectorized version of this is 
\[{vec(\hat{Y})}=(K\otimes A)[K\otimes A+\lambda I]^{-1}vec(Y)\]
 the matrix $A(\lambda):=(K\otimes A)[K\otimes A+\lambda I]^{-1}$ is the influence matrix used in the calculation of GCV.

\end{document}